\documentclass{amsart}
\usepackage{amssymb, latexsym}
\usepackage{amscd}

\newtheorem{remark}{Remark}
\newtheorem{proposition}{Proposition}
\newtheorem{definition}{Definition}
\newtheorem{example}{Example}

\begin{document}

\title{Local Tameness and Ends of Groups}
\author{Rita Gitik}
\email{ritagtk@umich.edu}
\address{Department of Mathematics \\ University of Michigan \\ Ann Arbor, MI, 48109}  

\date{\today}

\begin{abstract} 
We investigate connections between  local tameness of a group and a number of its ends.
\end{abstract}

\maketitle


\section{Introduction}

Let $H$ be a subgroup of a group $G$ given by a  presentation $G= \langle X|R \rangle $. 
Let $K$ be the standard presentation $2$-complex of $G$, i.e. $K$ has one vertex, $K$ has an edge, which is a loop, for every generator $x \in X$, and $K$ has a $2$-cell for every relator $r \in R$. The Cayley complex of $G$, denoted by $Cayley_2(G)$, is  the universal cover of $K$. Denote by $Cayley_2(G,H)$ the cover of $K$ corresponding to
a subgroup $H$ of $G$.

\begin{definition} cf. \cite{Gi1} and \cite{Mi}.
A finitely generated subgroup $H$ of a finitely presented group  $G$ is tame in $G$ if for any finite subcomplex $C$ of $Cayley_2(G,H)$ and for any component $C_0$ of $Cayley_2(G,H)-C$ the group $\pi_1(C_0)$ is finitely generated.
\end{definition}

Note that Definition 1 makes sense for any subgroup $H$ and for any presentation of a group $G$, however it was shown in \cite{Mi} that tameness of a finitely generated subgroup $H$ is independent of a finite presentation of a group $G$. As we are interested in investigating a property of a group, rather than a property of a presentation, we will work with finitely presented groups.

It is not known if there exists a finitely generated subgroup $H$ of a finitely presented group $G$ such that $H$ is not tame in $G$. Moreover, it is not known if the trivial subgroup is tame in any finite presentation of a group. However, an infinitely generated subgroup might not be tame in a finitely presented group.

\begin{example}
Let $F=\langle x,y \rangle$ be a free group of rank two, and let $F'$ be its commutator subgroup. The complex $Cayley_2(F,F')$ is one-dimensional and is homeomorphic to a lattice in the Euclidean plane, consisting of the horizontal lines $\{ y=n, n \in \mathbf{Z} \}$ and the vertical lines $\{ x=n, n \in \mathbf{Z} \}$. As the fundamental group of the complement of any finite connected subcomplex $C$ in $Cayley(F,F')$ is infinitely generated, $F'$ is not tame in $F$.
\end{example}

The following definition was given in \cite{Gi2}.

\begin{definition}
A  group $G$ is locally tame if all finitely generated subgroups of $G$ are tame in $G$.
\end{definition}

A finitely generated subgroup can be tame in an infinite presentation of a group. Moreover, an infinite presentation of a group can be locally tame.

\begin{example}
Let $F= \langle X \rangle $ be an infinitely generated free group and let $H$ be a finitely generated subgroup of $F$. Note that $Cayley_2(F,H)$ is one-dimensional. It has a finite subcomplex, called the core, which carries $H$, and the complement of the core in $Cayley_2(F,H)$ is a forest. Hence for any finite subcomplex of $Cayley_2(F,H)$ the fundamental group of any component of $Cayley_2(F,H) - C$ is finitely generated. Therefore any free group in standard free presentation is locally tame.
\end{example}

In this paper we discuss connections between local tameness and a number of ends of a group.

It was shown in \cite{Fr} and \cite{Ho} that any finitely generated group  has either $0,  1, 2$, or infinitely many ends.

\begin{remark}
A finitely presented  group has $0$ ends if and only if it is finite, so it is locally tame because its Cayley complex is finite.
\end{remark}

\begin{remark}
A finitely generated group  has two ends if and only if it is virtually infinite cyclic, i.e. it contains an infinite cyclic subgroup of finite index. 

It was shown in \cite{Gi2}  that virtually locally tame groups are locally tame and an infinite cyclic group is locally tame. Hence any finitely presented group with $2$ ends is locally tame. 
\end{remark}

The main result of this paper is the following theorem. The proof is given in Section 2.

\bigskip

\textbf{Theorem.}

\emph{If all finitely presented groups with one end are locally tame then all finitely presented groups  with infinitely many ends are locally tame.}

\begin{remark}
It was shown in \cite{Co} that a free product of two finitely generated groups has infinitely many ends unless both groups have order $2$, in which case it has $2$ ends.
\end{remark}

\bigskip

\section{Proof of the Theorem}

It was shown in \cite{Sta} that a finitely generated  group with infinitely many ends is either an amalgamated free product of finitely generated groups with finite amalgamating subgroup or an HNN-extension of a finitely generated group over a finite subgroup.

Dunwoody proved in \cite{Du} that a finitely presented group with infinitely many ends is a fundamental group of a finite graph of groups with all edge groups finite and all vertex groups finitely generated with at most one end. He also proved (private communication) that a finitely presented group with infinitely many ends is a fundamental group of a finite graph of groups with all edge groups finite and all vertex groups finitely presented with at most one end. The proof of the stronger result, apparently, did not appear in print.

\bigskip

The following proposition follows from the definition of a graph of groups.

\begin{proposition}
Let $G$ be the fundamental group of a finite graph of groups with all edge groups finitely generated.

If all the vertex groups of $G$ are finitely generated then the presentation of $G$ given by that graph of groups structure is finitely generated.

If all the vertex groups of $G$ are finitely presented then the presentation of $G$ given by that graph of groups structure is finite.
\end{proposition}

Recall that a group is Noetherian if all its subgroups are finitely generated.

The following result was proved by the author in \cite{Gi3}.

\begin{proposition}
Let $G$ be the fundamental group of a finite graph of groups with Noetherian edge groups such that the presentation of $G$ given by the graph of groups structure is finite. If all the vertex groups of $G$ are locally tame then $G$ is locally tame.
\end{proposition}

\textbf{Proof of the Theorem.}

Let $G$ be the fundamental group of a finite graph of groups with finite edge groups and finitely presented locally tame vertex groups.

 Proposition 1 implies that the presentation of $G$ given by that graph of groups structure is finite. As finite groups are Noetherian, Proposition 2 implies that $G$  is locally tame. 

The accessibility theorem of Dunwoody, quoted above, states that a finitely presented group with infinitely many ends is a fundamental group of a finite graph of groups with all edge groups finite and all vertex groups finitely presented with at most one end.

Remark 1 states that all finitely presented groups with $0$ ends are locally tame. Therefore, if all finitely  presented groups with one end are locally tame then all finitely presented groups with infinitely many ends are locally tame.  

\section{Open Questions}

\begin{remark}
The group $F_2 \times F_2$ is a source of counterexamples. Is it locally tame? Lemma $7$ in \cite{Gi2} implies that the trivial subgroup is tame in $F_2 \times F_2$. It is proven in \cite{Co} that the direct product of two finitely generated groups has one end, so $F_2 \times F_2$ has one end.
\end{remark}

Recall that a group has f.g.i.p. (finitely generated intersection property) if the intersection of any two of its finitely generated subgroups is finitely generated.

\begin{remark} Is the group $F_2 \times F_1$ locally tame? It is shown in \cite{Mo} that $F_2 \times F_1$ does not have f.g.i.p.. Indeed, let $F_2$ be generated by the elements $x$ and $y$ and let $F_1$ be generated by the element $t$. Let $K$ be a subgroup of $F_2 \times F_1$ generated by the elements $x$ and $ty$. Then $F_2 \cap K $ is the normal closure of $x$ in $F_2$, hence it is not finitely generated.
\end{remark}

\begin{remark} Let $A$ and $B$ be finitely presented locally tame groups. Is it true that a finitely generated subgroup of $A \times B$ is tame in $A \times B$  if its intersections with $A$ and $B$  are finitely generated.
\end{remark}

\begin{remark}
A positive answer to the question stated in the previous remark will imply that
the direct product of two finitely presented locally tame groups is locally tame if any finitely generated subgroup intersects each factor in a finitely generated subgroup.
\end{remark}

\begin{remark}
Let a finitely presented group $G$ be a direct product $G=N \times L$, where $L$ is finite and $N$ is finitely presented. Remark 10 implies that $G$ is locally tame if and only if $N$ is locally tame.
\end{remark}

\begin{remark}
The following result was proved in \cite{Gi2}. Let $K_0$ be a finite index subgroup of a finitely presented group $K$. A finitely generated subgroup $H$ of $K$ is tame in $K$ if and only if $H \cap K_0$ is tame in $K_0$. It follows that virtually locally tame groups are locally tame and finite index subgroups of locally tame groups are locally tame.
\end{remark}

\section{Acknowledgment}

The author would like to thank Martin Dunwoody, Ilya Kapovich, and Mike Mihalik for helpful comments.

\end{document}